\documentclass[12pt,twoside]{article}
\usepackage[english]{babel}
\usepackage[latin1]{inputenc}
\usepackage{amsmath}
\usepackage{amssymb,amsfonts}
\usepackage{graphicx}  %% for include pictures *.eps

\newcommand{\bR}{\mathbf{R}}

\newcommand{\bS}{\mathbf{S}}

\newcommand{\HYP}{\mathbb{H}^3}
\newcommand{\HYN}{\mathbb{H}^n}
\newcommand{\SXR}{\bS^2\!\times\!\bR}

\begin{document}
\pagestyle{myheadings}
\markboth{\centerline{Jen\H o Szirmai}}
{Decomposition method related to saturated hyperball packings}
\title
{Decomposition method related to saturated hyperball packings}

\author{\normalsize{Jen\H o Szirmai} \\
\normalsize Budapest University of Technology and \\
\normalsize Economics Institute of Mathematics, \\
\normalsize Department of Geometry \\
\date{\normalsize{\today}}}

\maketitle

%%%%%%%%%%%%%%%%%%%%%%%%%%%%%%%%%%%%%%%%%%%%

\begin{abstract}

In this paper we study the problem of hyperball (hypersphere) packings in 
$3$-dimensional hyperbolic space. We introduce a new definition of the non-compact saturated ball packings and describe to each saturated hyperball packing,
a new procedure to get a decomposition of 3-dimensional hyperbolic space $\HYP$ into truncated
tetrahedra. Therefore, in order to get a density upper bound for hyperball packings, it is sufficient to determine
the density upper bound of hyperball packings in truncated simplices.
\end{abstract}

%%%%%%%%%%%%%%%%%%%%%%%%%%%%%%%%%%%%%%%%%%%
\newtheorem{theorem}{Theorem}[section]
\newtheorem{corollary}[theorem]{Corollary}
\newtheorem{conjecture}{Conjecture}[section]
\newtheorem{lemma}[theorem]{Lemma}
\newtheorem{exmple}[theorem]{Example}
\newtheorem{defn}[theorem]{Definition}
\newtheorem{rmrk}[theorem]{Remark}
            %%% for no-italic, numbered environments, use:
\newenvironment{definition}{\begin{defn}\normalfont}{\end{defn}}
\newenvironment{remark}{\begin{rmrk}\normalfont}{\end{rmrk}}
\newenvironment{example}{\begin{exmple}\normalfont}{\end{exmple}}
            %%% for unnumbered environments, use f.e.
\newenvironment{acknowledgement}{Acknowledgement}
%%%%%%%%%%%%%%%%%%%%%%%%%%%%%%%%%%%%%%%%%%%%%%%%%%%%%%%%%%%%%%%%%%%%

%============================================================================%
%                             the main article                               %
%============================================================================%

%%%%%%%%%%%%%%%%%%%%%%%%%%%%%%%%%%%%%%%%%%%%%%%%%%%%%%%%%%%%%%%%%%%%%%%%%%%%%%
\section{Introduction}
In $n$-dimensional hyperbolic space $\HYN$ ($n \ge 2$) there are $3$ kinds
of generalized "balls (spheres)":  the usual balls (spheres), horoballs (horospheres) and hyperballs (hyperspheres).

The classical problems of ball packings and coverings with congruent {\it generalized balls} of hyperbolic spaces $\HYN$ are extensively discussed in the literature,
however there are several essential open questions e.g.:
\begin{enumerate}
\item What are the optimal ball packing and covering configurations of {\it usual spheres} and what are their densities ($n \ge 3$) (see \cite{Be,G--K--K, K98, MSz17})?
\item The monotonicity of the density related to the Böröczky type ball configurations depending on the radius of the congruent balls $(n \ge 4)$ (see \cite{B--F64, M}).
\item What are the optimal horoball packing and covering configurations and what are their densities allowing horoballs in different types ($n \ge 4$) (see \cite{B78, KSz, KSz14})?
\item What are the optimal packing and covering arrangements using non-compact balls (horoballs and hyperballs) and what are their densities?
These are the so-called hyp-hor packings and coverings (see \cite{Sz17-1}).
\item What are the optimal hyperball packing and covering configurations and what are their densities ($ n\ge 3$)?
\end{enumerate}
In this paper we study the $5^{th}$ question related to saturated, congruent hyperball packings in $3$-dimensional hyperbolic space $\HYP$.

In the hyperbolic plane $\mathbb{H}^2$ the universal upper bound of the hypercycle packing density is $\frac{3}{\pi}$,
proved by I.~Vermes in \cite{V79} and the universal lower bound of the hypercycle covering density is $\frac{\sqrt{12}}{\pi}$
determined by I.~Vermes in \cite{V81}.

In \cite{Sz06-1} and \cite{Sz06-2} we studied the regular prism tilings (simple truncated Coxeter orthoscheme tilings) and the corresponding optimal hyperball packings in
$\mathbb{H}^n$ $(n=3,4)$ and we extended the method developed in the former paper \cite{Sz13-3} to
5-dimensional hyperbolic space.
Moreover, their metric data and their densities have been determined. In paper \cite{Sz13-4} we studied the $n$-dimensional hyperbolic regular prism honeycombs
and the corresponding coverings by congruent hyperballs and we determined their least dense covering densities.
Furthermore, we formulated conjectures for the candidates of the least dense hyperball
covering by congruent hyperballs in the 3- and 5-dimensional hyperbolic space ($n \in \mathbb{N},3 \le n \le 5)$.

In \cite{Sz17-2} we discussed congruent and non-congruent hyperball (hypersphere) packings of the truncated regular tetrahedron tilings.
These are derived from the Coxeter simplex tilings $\{p,3,3\}$ $(7\le p \in \mathbb{N})$ and $\{5,3,3,3,3\}$
in $3$- and $5$-dimensional hyperbolic space.
We determined the densest hyperball packing arrangement and its density
with congruent hyperballs in $\mathbb{H}^5$ and determined the smallest density upper bounds of
non-congruent hyperball packings generated by the above tilings in $\HYN,~ (n=3,5)$.

In \cite{Sz17-1} we deal with the packings derived by horo- and hyperballs (briefly hyp-hor packings) in $n$-dimensional hyperbolic spaces $\HYN$
($n=2,3$) which form a new class of the classical packing problems.
We constructed in the $2-$ and $3-$dimensional hyperbolic spaces hyp-hor packings that
are generated by complete Coxeter tilings of degree $1$ i.e. the fundamental domains of these tilings are simple frustum orthoschemes
and we determined their densest packing configurations and their densities.
We proved using also numerical approximation methods that in the hyperbolic plane ($n=2$) the density of the above hyp-hor packings arbitrarily approximate
the universal upper bound of the hypercycle or horocycle packing density $\frac{3}{\pi}$ and
in $\HYP$ the optimal configuration belongs to the $\{7,3,6\}$ Coxeter tiling with density $\approx 0.83267$.
Furthermore, we analyzed the hyp-hor packings in
truncated orthosche\-mes $\{p,3,6\}$ $(6< p < 7, ~ p\in \mathbb{R})$ whose
density function is attained its maximum for a parameter which lies in the interval $[6.05,6.06]$
and the densities for parameters lying in this interval are larger that $\approx 0.85397$. That means that these
locally optimal hyp-hor configurations provide larger densities that the B\"or\"oczky-Florian density upper bound
$(\approx 0.85328)$ for ball and horoball packings but these hyp-hor packing configurations can not be extended to the entirety of hyperbolic space $\mathbb{H}^3$.

In \cite{Sz17} we studied a large class of hyperball packings in $\HYP$
that can be derived from truncated tetrahedron tilings.
In order to get a density upper bound for the above hyperball packings, it is sufficient
to determine this density upper bound locally, e.g. in truncated tetrahedra.
Thus we proved that if the truncated tetrahedron is regular, then the density
of the densest packing is $\approx 0.86338$. This is larger than the B\"or\"oczky-Florian density upper bound for balls and horoballs
but our locally optimal hyperball packing configuration cannot be extended to the entirety of
$\mathbb{H}^3$. However, we described a hyperball packing construction,
by the regular truncated tetrahedron tiling under the extended Coxeter group $\{3, 3, 7\}$ with maximal density $\approx 0.82251$.

Recently, (to the best of author's knowledge) the candidates for the densest hyperball
(hypersphere) packings in the $3,4$ and $5$-dimensional hyperbolic space $\mathbb{H}^n$ are derived by the regular prism
tilings which have been in papers \cite{Sz06-1}, \cite{Sz06-2} and \cite{Sz13-3}.

{\it In this paper we study hyperball (hypersphere) packings in
$3$-dimensional hyperbolic space.
We develop a decomposition algorithm that for each saturated hyperball packing provides a decomposition of $\HYP$
into truncated tetrahedra. Therefore, in order to get a density upper bound for hyperball packings, it is sufficient to determine
the density upper bound of hyperball packings in truncated simplices.}
\section{The projective model and saturated hyperball packings of $\HYP$}
We use for $\mathbb{H}^3$ (and analogously for $\HYN$, $n\ge3$) the projective model in the Lorentz space $\mathbb{E}^{1,3}$
that denotes the real vector space $\mathbf{V}^{4}$ equipped with the bilinear
form of signature $(1,3)$,
$
\langle \mathbf{x},~\mathbf{y} \rangle = -x^0y^0+x^1y^1+x^2y^2+ x^3 y^3,
$
where the non-zero vectors
$
\mathbf{x}=(x^0,x^1,x^2,x^3)\in\mathbf{V}^{4} \ \  \text{and} \ \ \mathbf{y}=(y^0,y^1,y^2,y^3)\in\mathbf{V}^{4},
$
are determined up to real factors, for representing points of $\mathcal{P}^n(\mathbb{R})$. Then $\mathbb{H}^3$ can be interpreted
as the interior of the quadric
$
Q=\{(\mathbf{x})\in\mathcal{P}^3 | \langle  \mathbf{x},~\mathbf{x} \rangle =0 \}=:\partial \mathbb{H}^3
$
in the real projective space $\mathcal{P}^3(\mathbf{V}^{4},
\mbox{\boldmath$V$}\!_{4})$ (here $\mbox{\boldmath$V$}\!_{4}$ is the dual space of $\mathbf{V}^{4}$).
Namely, for an interior point $\mathbf{y}$ holds $\langle  \mathbf{y},~\mathbf{y} \rangle <0$.

Points of the boundary $\partial \mathbb{H}^3 $ in $\mathcal{P}^3$
are called points at infinity, or at the absolute of $\mathbb{H}^3 $. Points lying outside $\partial \mathbb{H}^3 $
are said to be outer points of $\mathbb{H}^3 $ relative to $Q$. Let $(\mathbf{x}) \in \mathcal{P}^3$, a point
$(\mathbf{y}) \in \mathcal{P}^3$ is said to be conjugate to $(\mathbf{x})$ relative to $Q$ if
$\langle \mathbf{x},~\mathbf{y} \rangle =0$ holds. The set of all points which are conjugate to $(\mathbf{x})$
form a projective (polar) hyperplane
$
pol(\mathbf{x}):=\{(\mathbf{y})\in\mathcal{P}^3 | \langle \mathbf{x},~\mathbf{y} \rangle =0 \}.
$
Thus the quadric $Q$ induces a bijection
(linear polarity $\mathbf{V}^{4} \rightarrow
\mbox{\boldmath$V$}\!_{4})$
from the points of $\mathcal{P}^3$ onto their polar hyperplanes.

Point $X (\bold{x})$ and hyperplane $\alpha (\mbox{\boldmath$a$})$
are incident if $\bold{x}\mbox{\boldmath$a$}=0$ ($\bold{x} \in \bold{V}^{4} \setminus \{\mathbf{0}\}, \ \mbox{\boldmath$a$}
\in \mbox{\boldmath$V$}_{4}
\setminus \{\mbox{\boldmath$0$}\}$).

The hypersphere (or equidistance surface) is a quadratic surface at a constant distance
from a plane (base plane) in both halfspaces. The infinite body of the hypersphere, containing the base plane, is called hyperball.

The {\it half hyperball } with distance $h$ to a base plane $\beta$ is denoted by $\mathcal{H}^h_+$.
The volume of a bounded hyperball piece $\mathcal{H}^h_+(\mathcal{A})$,
delimited by a $2$-polygon $\mathcal{A} \subset \beta$, and its prism orthogonal to $\beta$, can be determined by the classical formula
(2.1) of J.~Bolyai \cite{B31}.
\begin{equation}
\mathrm{Vol}(\mathcal{H}^h_+(\mathcal{A}))=\frac{1}{4}\mathrm{Area}(\mathcal{A})\left[k \sinh \frac{2h}{k}+
2 h \right], \tag{2.1}
\end{equation}
The constant $k =\sqrt{\frac{-1}{K}}$ is the natural length unit in
$\mathbb{H}^3$, where $K$ denotes the constant negative sectional curvature. In the following we may assume that $k=1$.

Let $\mathcal{B}^h$ be a hyperball packing in $\HYP$ with congruent hyperballs of height $h$.

The notion of {\it saturated packing} follows from that fact that the density of any packing can be improved by adding further packing elements as long as there is
sufficient room to do so. However, we usually apply this notion for packings with congruent elements.
\begin{figure}[ht]
\centering
\includegraphics[width=5.5cm]{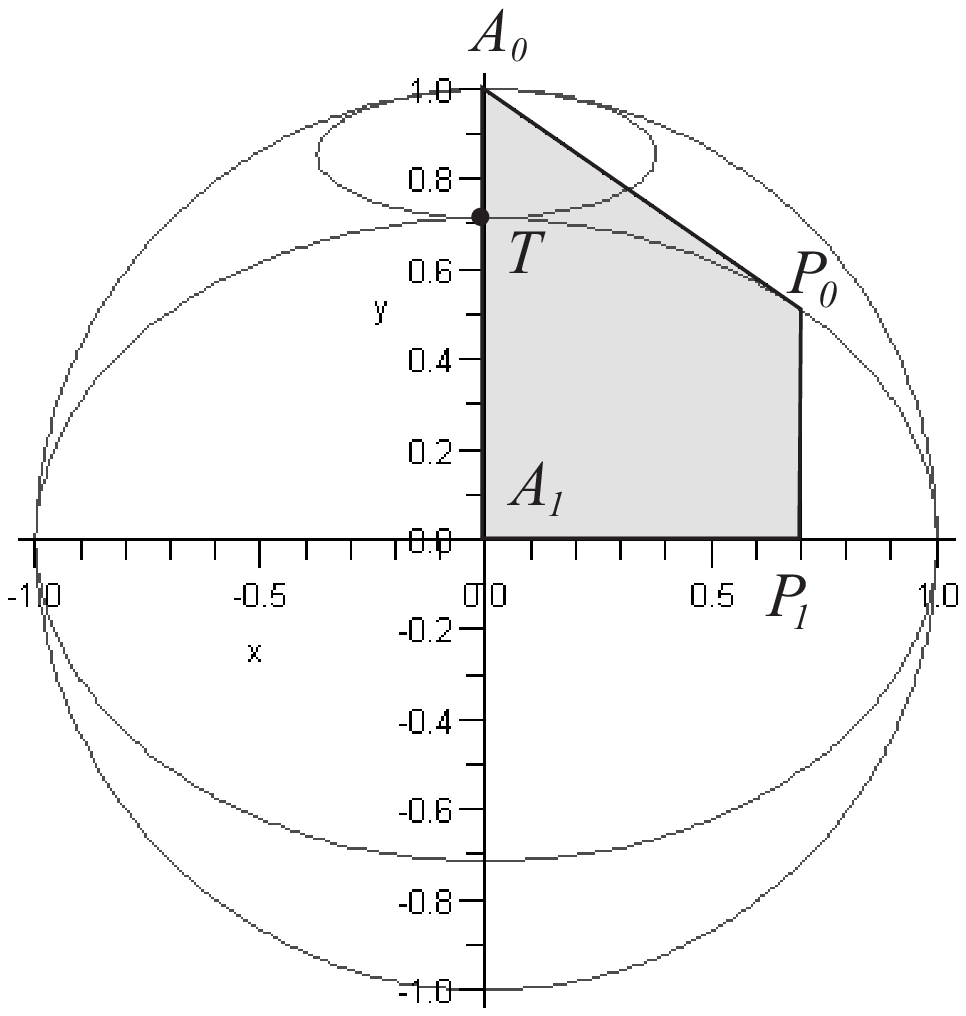}~ ~ \includegraphics[width=5.5cm]{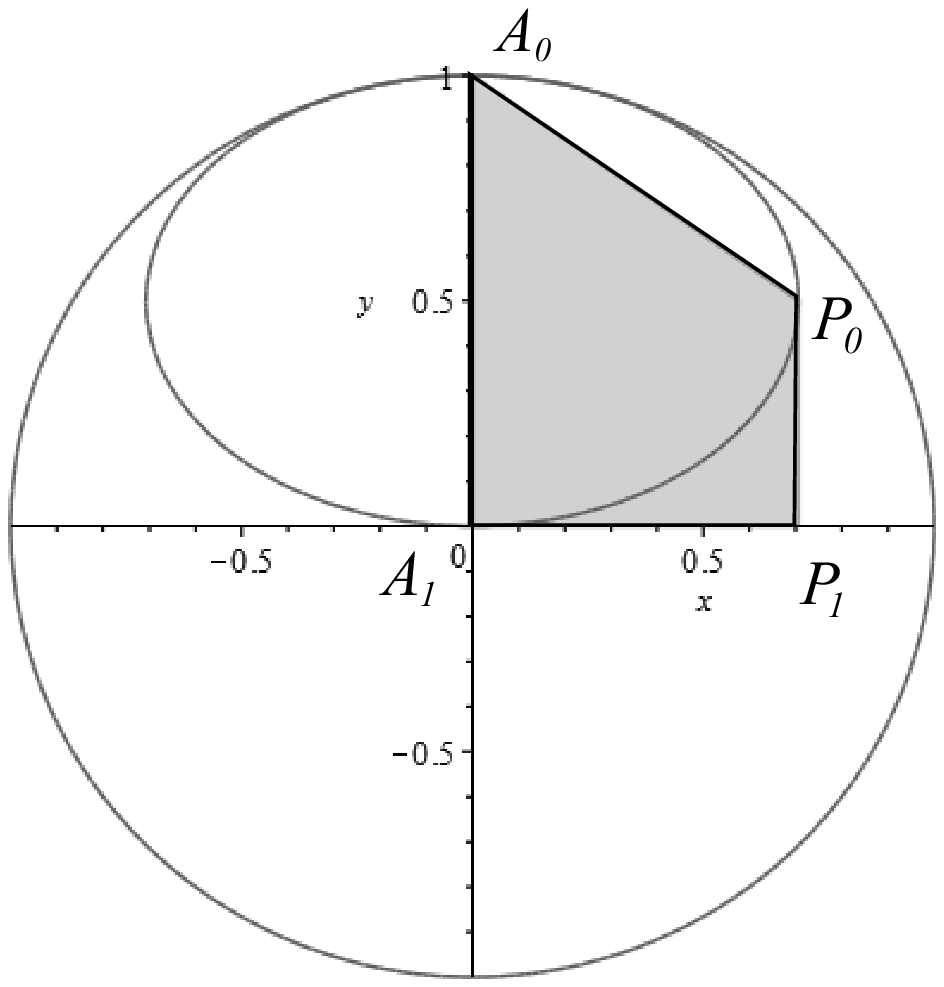}
\caption{a. Saturated hyp-hor packing, at present $a=0.7$.~ b. Saturated horocycle packing with parameter $a=\frac{1}{\sqrt{2}}$.}
\label{}
\end{figure}
Now, {\it we modify the classical definition of saturated packing} for non-compact ball packings in $n$-dimensional hyperbolic space $\HYN$ ($n\ge 2$ integer parameter):
\begin{defn}
A ball packing with non-compact balls (horoballs or/and hyperballs) in $\HYN$ is saturated if no new non-compact ball can be added to it.
\end{defn}
We illustrate the meaning of the above definition by $2$-dimensional Coxeter tilings given by the Coxeter symbol $[\infty]$ (see Fig.~1),
which are denoted by $\mathcal{T}_a$. The fundamental domain of $\mathcal{T}_a$ is a Lambert quadrilateral $A_0A_1P_0P_1$ (see  \cite{Sz17-1})
that is denoted by $\mathcal{F}_a$. It is derived by the truncation of the orthoscheme $A_0A_1A_2$
by the polar line $\pi$ of the outer vertex $A_2$. The other initial principal vertex $A_0$ of the orthoscheme is lying on the absolute quadric of the Beltrami-Cayley-Klein model

The images of $\mathcal{F}_a$ under reflections on its sides fill the hyperbolic plane $\mathbb{H}^2$ without overlap. The tilings $\mathcal{T}_a$ contain a free parameter
$0 < a < 1,~ ~ a \in \mathbb{R}$. The
polar straight line of $A_2$ is $\pi$ and $\pi \cap A_0A_2=P_0$, $\pi \cap A_1A_2=P_1$. If we fix the parameter $a$ then a optimal hypercycle tiling can be derived from
the mentioned Coxeter tiling (see Fig.~1.a) but here there are sufficient rooms to add horocycles with centre $A_0$ and with centres at the images of $A_0$.
This saturated {\it hyp-hor packing} (packing with horo- and hyperballs) is illustrated in Fig.~1.a. The Fig.~1.b shows a saturated horocycle packing belonging to the same Coxeter tiling.

To obtain hyperball (hypersphere) packing bounds it obviously suffices to study saturated hyperball packings (using the above definition) and in what follows we assume that
all packings are saturated unless otherwise stated.
\section{Decomposition into truncated tetrahedra}
We take the set of hyperballs $\{ \mathcal{H}^h_i\}$ of a saturated hyperball packing $\mathcal{B}^h$ (see Definition 2.1).
Their base planes are denoted by $\beta_i$.
Thus in a saturated hyperball packing the distance between two ultraparallel base planes
$d(\beta_i,\beta_j)$ is at least $2h$ (where for the natural indices holds $i < j$
and $d$ is the hyperbolic distance function).

In this section we describe a procedure to get a decomposition of 3-dimensional hyperbolic space $\HYP$ into truncated
tetrahedra corresponding to a given saturated hyperball packing.
\begin{enumerate}
\item The notion of the radical plane (or power plane) of two Euclidean spheres can be extended to the hyperspheres.
The radical plane (or power plane) of two non-intersecting hyperspheres
is the locus of points at which tangents drawn to both hyperspheres have the same length (so these points have equal power with respect
to the two non-intersecting hyperspheres). If the two non-intersecting hyperspheres are congruent also in Euclidean sense in the model then their 
radical plane coincides with their "Euclidean symmetry plane" and any two congruent hypersphere can be transformed into such an hypersphere arrangement.

Using the radical planes of the hyperballs $\mathcal{H}^h_i$, similarly to the Euclidean space, can be constructed
the unique Dirichlet-Voronoi (in short $D-V$) decomposition of $\HYP$ to the given hyperball packing $\mathcal{B}^h$. Now, the
$D-V$ cells are infinite hyperbolic polyhedra containing the corresponding hyperball,
and its vertices are proper points of $\mathbb{H}^3$. We note here (it is easy to see), that a vertex of any $D-V$ cell
cannot be outer or boundary point of $\mathbb{H}^3$ relative to $Q$, because the hyperball packing $\mathcal{B}^h$ is saturated by the Definition 2.1.

\item We consider an arbitrary {\it proper} vertex $P \in \HYP$ of the above $D-V$ decomposition and the hyperballs $\mathcal{H}^h_i(P)$
whose $D-V$ cells
meet at $P$. The base planes of the hyperballs $\mathcal{H}^h_i(P)$ are denoted by $\beta_i(P)$, and these planes determine a
non-compact polyhedron $\mathcal{D}^i(P)$ with the intersection of their halfspaces
containing the vertex $P$. Moreover, denote $A_1,A_2,A_3,\dots$ the outer vertices of $\mathcal{D}^i(P)$ and cut off
$\mathcal{D}^i(P)$ with the polar planes $\alpha_j(P)$ of its outer vertices $A_j$. Thus, we obtain a convex compact polyhedron
$\mathcal{D}(P)$.
This is bounded by the base planes $\beta_i(P)$ and "polar planes" $\alpha_i(P)$. Applying this procedure for all vertices of the
above Dirichlet-Voronoi decomposition, we obtain an other decomposition of $\HYP$ into convex polyhedra.
\item We consider $\mathcal{D}(P)$ as a tile of the above decomposition. The planes from the finite set of base planes $\{\beta_i(P)\}$ are
called adjacent if there is a vertex $A_s$ of $\mathcal{D}^i(P)$ that lies on each of the above plane.
We consider non-adjacent planes $\beta_{k_1}(P),\beta_{k_{2}}(P),\beta_{k_{3}}(P), \dots \beta_{k_{m}}(P) \in \{\beta_i(P)\}$ $(k_l  \in \mathbb{N}^+, ~ l=1,2,3,\dots m)$
that have an outer point of intersection denoted by $A_{k_1\dots k_m}$. Let $N_{\mathcal{D}(P)}$ denote the {\it finite} number of the outer points $A_{k_1\dots k_m}$
related to $\mathcal{D}(P)$. It is clear, that its
minimum is $0$ if $\mathcal{D}^i(P)$ is tetrahedron.
The polar plane $\alpha_{k_1\dots k_m}$ of $A_{k_1\dots k_m}$ is orthogonal to planes $\beta_{k_1}(P),\beta_{k_2}(P), \dots \beta_{k_m}(P)$ 
(thus it contain their poles $B_{k_1}$, $B_{k_2}$,\dots $B_{k_m}$) and divides $\mathcal{D}(P)$ into two convex polyhedra
$\mathcal{D}_1(P)$ and $\mathcal{D}_2(P)$.
\item If $N_{\mathcal{D}_1(P)} \ne 0$ and $N_{\mathcal{D}_2(P)} \ne 0$ then $N_{\mathcal{D}_1(P)} < N_{\mathcal{D}(P)}$ and $N_{\mathcal{D}_2(P)} < N_{\mathcal{D}(P)}$ then
we apply the point 3 for polyhedra $\mathcal{D}_i(P),~ i \in \{1,2\}$.
\item If $N_{\mathcal{D}_i(P)} \ne 0$ or $N_{\mathcal{D}_j(P)} = 0$ ($i \ne j,~ i,j\in\{1,2\}$) then we consider the polyhedron $\mathcal{D}_i(P)$ where
$N_{\mathcal{D}_i(P)}=N_{\mathcal{D}(P)}-1$ because the vertex $A_{k_1\dots k_m}$ is left out and apply the point 3.
\item If $N_{\mathcal{D}_1(P)}=0$ and $N_{\mathcal{D}_2(P)}=0$ then the procedure is over for $\mathcal{D}(P)$. We continue the procedure with the next cell.
\item It is clear, that the above plane $\alpha_{k_1\dots k_m}$ intersects every hyperball $ \mathcal{H}^h_j(P)$, $(j=k_1\dots k_m)$.
\begin{lemma}
The plane $\alpha_{k_1 \dots k_m}$ of $A_{k_1 \dots k_m}$ does not intersect the hyperballs $\mathcal{H}^h_s(P)$ where $A_{k_1 \dots k_m} \notin \beta_s(P).$
\end{lemma}
{\bf Proof}

Let $\mathcal{H}^h_s(P)$, ($A_{k_1 \dots k_m} \notin \beta_s(P)$) be an arbitrary hyperball corresponding to $\mathcal{D}(P)$ with base plane $\beta_s(P)$ whose pole
is denoted by $B_s$.
The common perpendicular $\sigma$ of the planes $\alpha_{k_1 \dots k_m}$ and $\beta_s(P)$ is the line through the point $A_{k_1 \dots k_m}$ and $B_s$.
We take a plane $\kappa$ containing the above common perpendicular, and its intersections with $\mathcal{D}(P)$ and $\mathcal{H}^h_s(P)$ are denoted by
$\phi$ and $\eta$. We
obtain the arrangement illustrated in Fig.~2 which coincides with the solution investigated in \cite{V79}. There I.~Vermes noticed that
the straight line $\phi=\alpha_{k_1 \dots k_m}\cap \kappa$ does not intersect the hypercycle $\eta=\mathcal{H}^h_s(P)\cap \kappa$.
The plane $\alpha_{k_1 \dots k_m}$ and the hyperball $\mathcal{H}^h_s(P)$ can be generated by rotation of $\phi$ and $\eta$ about
the common perpendicular $\sigma$; therefore, they are disjoint. $\square$
\begin{figure}[ht]
\centering
\includegraphics[width=12cm]{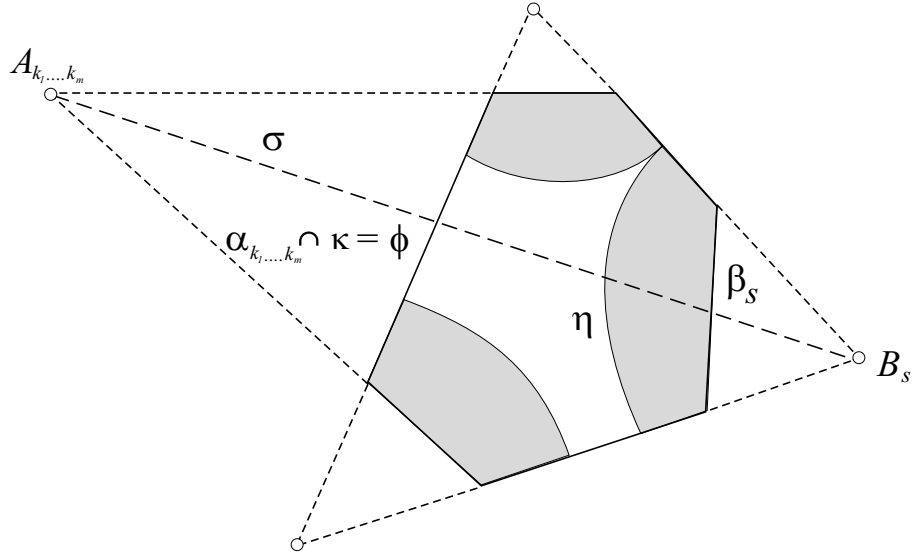}
\caption{The plane $\kappa$ and its intersections with $\mathcal{D}(P)$ and $\mathcal{H}^h_s(P)$}
\label{}
\end{figure}
\item We have seen in steps 3, 4, 5 and 6 that the number of the outer vertices $A_{k_1 \dots k_m}$ of any polyhedron obtained after the cutting process is less than the original one, and
we have proven in step 7 that the original hyperballs form packings in the new polyhedra $\mathcal{D}_1(P)$ and $\mathcal{D}_2(P)$, as well.
We continue the cutting procedure described in step 3 for both polyhedra $\mathcal{D}_1(P)$ and $\mathcal{D}_2(P)$. If a derived polyhedron is a truncated
tetrahedron then the cutting procedure does not give new polyhedra, thus the procedure will not be continued.
Finally, after a {\it finite number of cuttings} we get a decomposition of $\mathcal{D}(P)$ into truncated tetrahedra,
and in any truncated tetrahedron the corresponding congruent hyperballs from $\{ \mathcal{H}^h_i\}$ form a packing. Moreover, we apply the above method for the further cells.
\end{enumerate}
Finally we get the following
\begin{theorem}
The above described algorithm provides for each congruent saturated hyperball packing a decomposition of $\HYP$ into truncated tetrahedra. ~ ~$\square$
\end{theorem}
The above procedure is illustrated for regular octahedron tilings derived by the regular prism tilings with Coxeter-Schl\"afli symbol $\{p,3,4\}$, $6<p \in \mathbb{N}$.
These Coxeter tilings and
the corresponding hyperball packings are investigated in \cite{Sz06-1}.
In this situation the convex polyhedron $\mathcal{D}(P)$ is a truncated octahedron (see Fig.~3) whose vertices $B_i$, $(i=1,2,3,4,5,6)$ are outer points
and the octahedron is cut off with their polar planes $\beta_i$. These planes are the base planes of the hyperballs $\mathcal{H}^h_i$.
We can assume that the centre of the octahedron coincides with the centre of the model.

First, we choose three non-adjacent base planes $\beta_2,\beta_3, \beta_4$. Their common point, denoted by $A_{234}$ and its polar plane
$\alpha_{234}$ are
determined by points $B_2,B_3,B_4$ containing the centre $P$ as well. Then
we consider the non-adjacent base planes $\beta_2,\beta_4, \beta_5$ and
the polar plane $\alpha_{245}$ of their common point $A_{245}$. It is clear that the points $B_2,B_4,B_5$ lie in the plane $\alpha_{245}$ (see Fig.~3).

By the above two "cuttings" we get the decomposition of $\mathcal{D}(P)$ into truncated simplices.
\begin{figure}[ht]
\centering
\includegraphics[width=6.5cm]{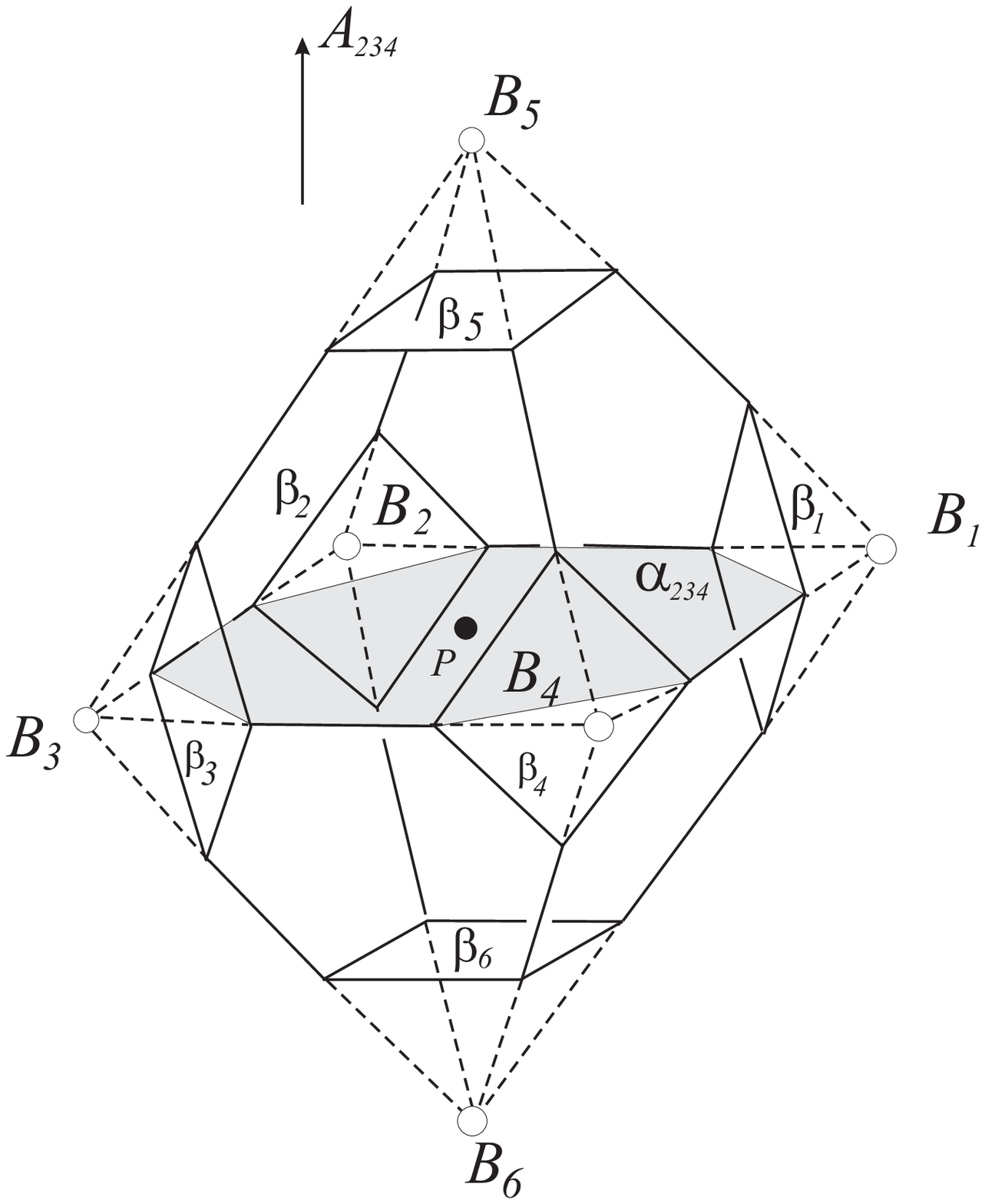} \includegraphics[width=6.5cm]{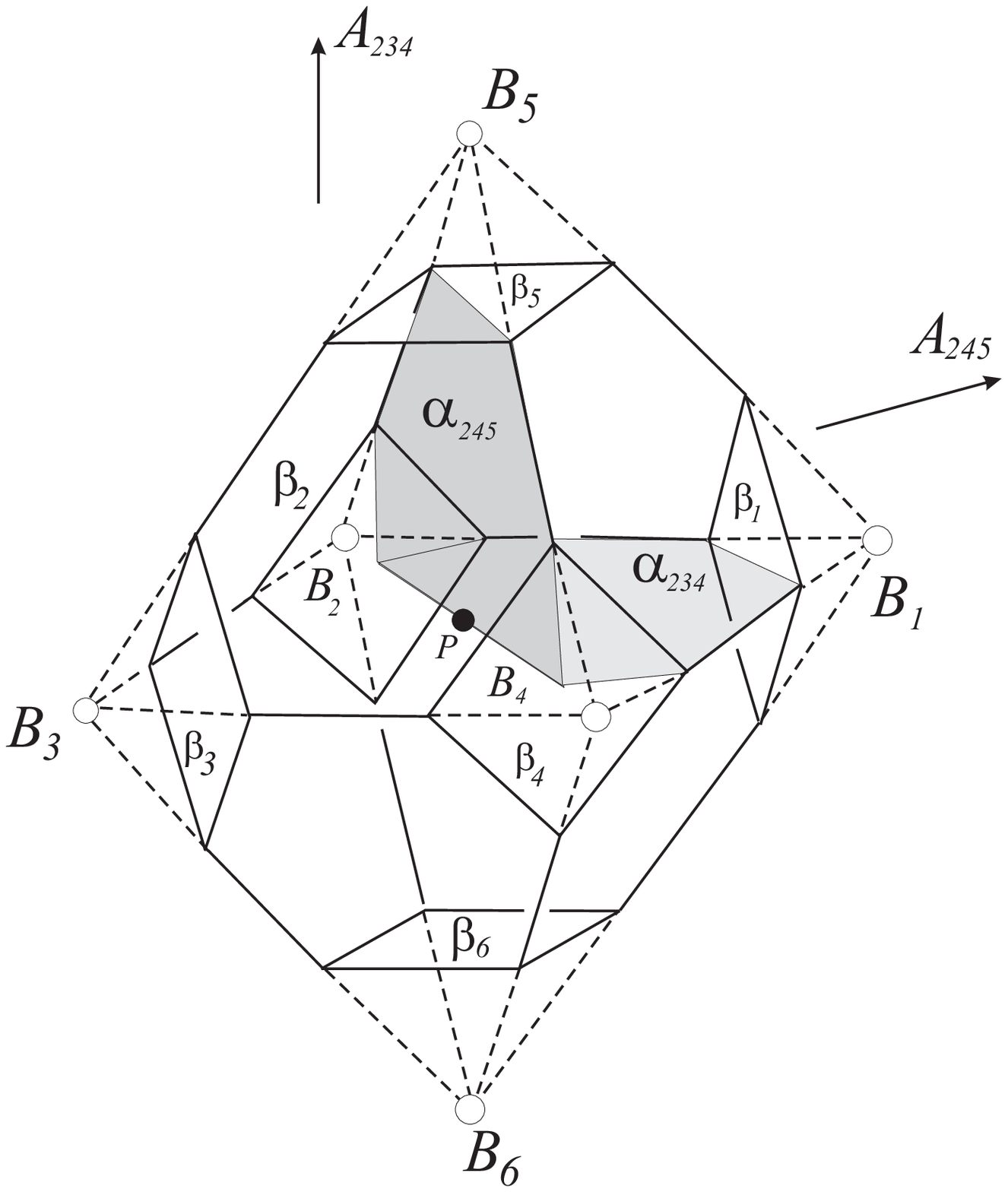}
\caption{Truncated octahedron tiling derived from the regular prism tilings with Coxeter-Schl\"afli symbol $\{p,3,4\}$ and its decomposition into
truncated tetrahedra}
\label{}
\end{figure}
\begin{rmrk}
\begin{enumerate}
\item
From the above section it follows that, to each saturated hyperball packing $\mathcal{B}^h$ of hyperballs $\mathcal{H}^h_i$ there is a
decomposition of $\HYP$ into truncated tetrahedra.
Therefore, in order to get a density upper bound for hyperball packings, it is sufficient to determine
the density upper bound of hyperball packings in truncated simplices.

We observed in \cite{Sz17} that some extremal properties of hyperball packings naturally belong to the regular truncated tetrahedron (or simplex, in general, see
Lemma 3.2 and Lemma 3.3 in \cite{Sz17}). Therefore, we studied hyperball
packings in regular truncated tetrahedra, and prove that if the truncated tetrahedron is regular, then the density of the densest packing is $\approx  0.86338$ (see
Theorem 5.1 in \cite{Sz17}). However, these hyperball packing configurations are only locally optimal, and cannot be extended to the whole space $\HYP$. Moreover, we
showed that the densest known hyperball packing, dually related to the regular prism tilings, introduced in \cite{Sz06-1}, can be realized by a regular truncated
tetrahedron tiling with density $\approx 0.82251$.
\item In \cite{Sz17-2} we discussed the problem of congruent and non-congruent hyperball
(hypersphere) packings to each truncated regular tetrahedron tiling.
These are derived from the Coxeter simplex tilings $\{p,3,3\}$ and $\{5,3,3,3,3\}$
in the 3 and 5-dimensional hyperbolic space. We determined the densest hyperball
packing arrangement and its density with congruent hyperballs in $\mathbb{H}^5$ $(\approx 0.50514)$
and determined the smallest density upper bounds of non-congruent hyperball
packings generated by the above tilings: in $\mathbb{H}^3$, $(\approx 0.82251$); in $\mathbb{H}^5$, $(\approx 0.50514$).
\end{enumerate}
\end{rmrk}

The question of finding the densest hyperball packings and horoball packings with horoballs
of different types in the $n$-dimensional hyperbolic spaces $n \ge 3$ has not been settled yet either (see e.g. \cite{Sz17}, \cite{KSz}, \cite{KSz14}).

Optimal sphere packings in other homogeneous Thurston geometries represent
another huge class of open mathematical problems. For these non-Euclidean geometries
only very few results are known (e.g. \cite{Sz11-2}, \cite{Sz14-1}). Detailed studies are the objective of ongoing research.
The applications of the above projective method seem to be interesting in (non-Euclidean)
crystallography as well, a topic of much current interest.

%============================================================================%
%                                references                                  %
%============================================================================%

\noindent
\footnotesize{Budapest University of Technology and Economics, Institute of Mathematics, \\
Department of Geometry, \\
H-1521 Budapest, Hungary. \\
E-mail:~szirmai@math.bme.hu \\
http://www.math.bme.hu/ $^\sim$szirmai}

\end{document}